\input epsf
\documentclass[a4paper,11pt]{amsart}
\usepackage{latexsym}
\newtheorem{thm}{Theorem}[section]
\newtheorem{lem}[thm]{Lemma}

\newtheorem{prop}[thm]{Proposition}

\newtheorem{rem}[thm]{Remark} 

\title{The Orchard Morphism}
\author{Roland Bacher \\
}

\begin{document}

\maketitle         
 
{\begin{abstract}
We define and prove uniqueness of a 
natural homomorphism (called the Orchard morphism)
from some groups associated naturally to a finite set $E$ to the
group ${\mathcal E}(E)$ of two-partitions of $E$ representing 
equivalence relations having at most two classes on $E$. 

As an application, we exhibit a natural equivalence relation on the set 
of points of generic finite configurations in ${\mathbf R}^d$.
\end{abstract} }


\section{Introduction}

Let $E$ be a set. The set ${\mathcal E}(E)$ of all equivalence relations on
$E$ into at most two equivalence classes can be endowed with a group structure,
which we call the {\it group of two-partitions} since its
elements represent partitions of $E$ into at most two disjoint subsets. 
Elements of ${\mathcal E}(E)$ can be represented by
functions $f:E\longrightarrow \{\pm 1\}$, well-defined up to multiplication
by $-1$. Such a function endows $E$ with the the equivalence relation 
given by the equivalence classes $f^{-1}(1)$ and $f^{-1}(-1)$. 

We denote by
$$E^{(l)}=\{(x_1,\dots,x_l)\in E^l\ \vert\ x_i\not=x_j \hbox{ for }  
1\leq i< j\leq l\}$$
the set of all sequences of length $l$ without repetitions with values in $E$.
Consider the multiplicative group $\{\pm 1\}^{E^{(l)}}$ formed
by all functions $\varphi: E^{(l)}\longrightarrow \{\pm 1\}$.
The symmetric group $\hbox{Sym}(\{1,\dots,l\})$ acts
on such functions by permutations of the $l$ arguments.
The set ${\mathcal F}_{+}(E^{(l)})$ of all symmetric
functions on which this action is trivial
is a subgroup of $\{\pm 1\}^{E^{(l)}}$ while the set ${\mathcal F}_{-}(E^{(l)})$
of all antisymmetric 
functions on which  $\hbox{Sym}(\{1,\dots,l\})$ acts by the 
the signature homomorphism is a free ${\mathcal F}_{+}(E^{(l)})-$set.
Since the product of two antisymmetric functions is always symmetric,
the union ${\mathcal F}_{\pm }(E^{(l)})={\mathcal F}_{+}(E^{(l)})\cup
{\mathcal F}_{-}(E^{(l)})$ is a subgroup of $\{\pm 1\}^{E^{(l)}}$.

Suppose now that the set $E$ is finite. The main 
result of this paper is the existence of a non-trivial
natural homomorphism
from the finite group ${\mathcal F}_{\pm }(E^{(l)})$
to the finite group ${\mathcal E}(E)$ of two-partitions
on $E$.
Naturality means that this homomorphism is 
$\hbox{Sym}(E)-$equivariant with respect to the obvious actions
by automorphisms of $\hbox{Sym}(E)$ on both groups 
${\mathcal F}_{\pm }(E^{(l)})$ and ${\mathcal E}(E)$.

We call this homomorphism the {\it Orchard homomorphism}. The Orchard 
homomorphism is the unique natural homomorphism
from ${\mathcal F}_{\pm }(E^{(l)})$ to ${\mathcal E}(E)$ which
is non-trivial for $1\leq l< \sharp(E)$. 
There exists however a natural non-trivial homomorphism
${\mathcal F}_{\pm}(E^{(2)})\longrightarrow{\mathcal E}(E)$
distinct from the Orchard homomorphism (which is trivial in this case)
if $E$ has $2$ elements. 

The existence of such a homomorphism from the set of symmetric
functions ${\mathcal F}_{+ }(E^{(l)})$ to ${\mathcal E}(E)$ is not
surprising. Its natural extension to the set
${\mathcal F}_{- }(E^{(l)})$ of antisymmetric functions is however
not completely obvious.

Natural examples of antisymmetric functions $E^{(l)}\longrightarrow
\{\pm 1\}$ where $E$ is a finite set arise for instance by
considering finite generic subsets of points in the oriented
real affine space ${\mathbf R}^{l-1}$ where generic means that any set
of $k\leq l$ points in $E$ is affinely independent. One gets
an antisymmetric function on $E^{(l)}$ by considering 
the orientation of simplices spanned by $l$ linearly ordered points of $E$.

The case $l=3$ may be illustrated as giving a natural rule to
plant trees of two distinct species in an orchard: The Queen of Heart
has randomly choosen $n$ generic locations $E$ in her future royal orchard and
asks Alice to plant cherry- and plum-trees in a natural way. Alice
assigns to each cyclically oriented triplet $(a,b,c)$
of points in $E$ the value
$1$ (respectively $-1$)
if the points $(a,b,c)$ define a positively (respectively negatively)
oriented triangle (where Wonderland is supposed to be an oriented plane).
This yields an antisymmetric function on $E^{(3)}$ whose image
by the Orchard morphism is a two-partition prescribing the choice
of the species (up to global permutation of all cherry- and plum-trees).

Figure 1 shows an example. The two-partition obtained by the Orchard morphism
(given by Proposition \ref{geomdescr} with $d=2$)
yields for the choosen nine positions three trees of one species
and six trees of the remaining species. The paper \cite{BG} contains 
many more examples, some of which are monochromatic (all vertices 
belong to the same common equivalence class).

\medskip

\centerline{\epsfysize8cm\epsfbox{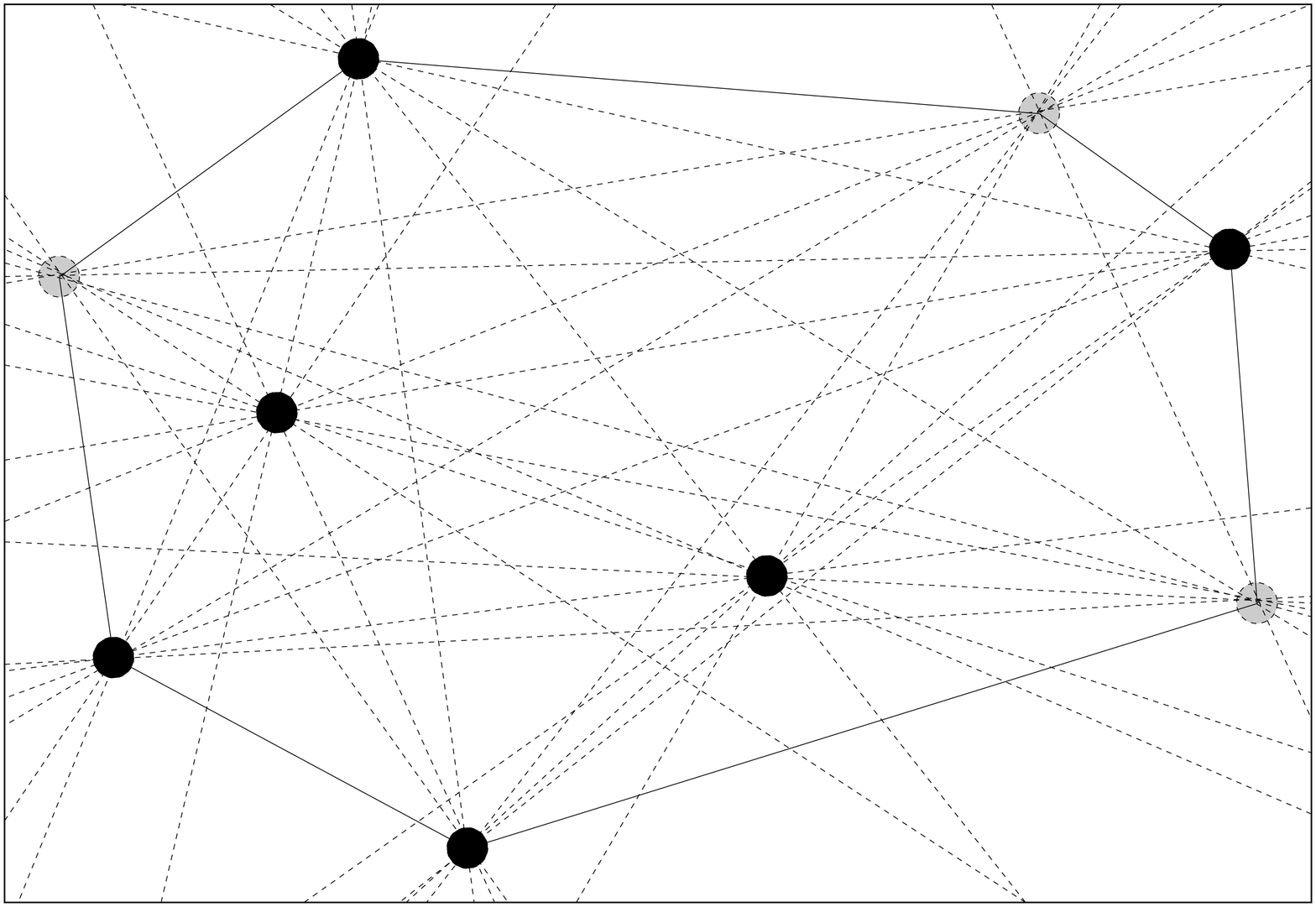}}
\centerline{Figure 1:
An orchard having 3 cherry- and 6 plumtrees in generic positions} 
\medskip

Finally, the Orchard morphism exists also in the case where the finite
set $E$ is endowed with a natural fixpoint-free involution 
$\iota:E\longrightarrow E$ which can
be thought of as a kind of orientation.
We call such a set orientable. Given an orientable set $(E,\iota)$ it
is natural to consider only structures on $E$
(equivalence relations, sets of functions etc.)
which are invariant (perhaps up to a sign) under the involution $\iota$.
We define in this setting analogues of the groups
${\mathcal E}(E)$ and ${\mathcal F}_\pm(E^{(l)})$ considered above and
construct the corresponding Orchard morphism.


\section{Two-partitions}

A {\it two-partition} is an unordered partition $\{A,B\}$ 
of a set $E=A\cup B$ into at most
two disjoint subsets. Two-partitions are the same as
equivalence relations having at most two classes. We 
will move freely between these two interpretations of two-partitions.
The word \lq\lq class'' will often be used instead of \lq\lq part
of the two-partition'' and a two-partition $\{A,B\}$ of $E$ will generally 
be written as $A\cup B$ or $E=A\cup B$.

A two-partition $E=A\cup B$
can be given by a pair $\pm \alpha$ of opposite functions where
$$\alpha:E\longrightarrow \{\pm 1\}$$ is defined by
$\alpha^{-1}(1)=A$ and $\alpha^{-1}(-1)=B$. The set ${\mathcal E}(E)$
of all such two-partitions is a vector space over
the field ${\mathbf F}_2$ of two elements. Its dimension is  
$\sharp(E)-1$ if $E$ is a finite set.
The pair $\pm 1$ of constant functions 
represents the identity and the group law $(\pm \alpha)(\pm \beta)$
is the obvious product $\pm \alpha\beta$ of functions. 
Set-theoretically, the product $(A_1\cup A_2)(B_1\cup B_2)$ 
of $2$ two-partitions on a set $E$ is given by $E=C_1\cup C_2$ where 
$C_1=\left(A_1\cap B_1\right)\cup\left(A_2\cap B_2\right)$ and
$C_2=\left(A_1\cap B_2\right)\cup\left(A_2\cap B_1\right)$.

Consider a simple graph $\Gamma$ (not necessarily finite) with 
vertices $V$ and unoriented edges $E$. 
Its adjacency matrix $A$ is the symmetric
matrix with rows and columns indexed by elements of $V$. All its entries are
zero except $A_{v,w}=1$ where $v\not= w$ are adjacent vertices of 
$\Gamma$ (i.e. $\{u,v\}$ is an edge of $\Gamma$). 

Our main tool in what follows is the following 
trivial and probably well-known observation:

\begin{lem} \label{equivgraphe}
Let $\Gamma$ be a simple graph with adjacency matrix $A$.
Suppose that there exists a constant $\gamma\in\{0,1\}$ such that 
$$A_{u,v}+A_{v,w}+A_{u,w}\equiv \gamma \pmod 2$$
for all triplets $(u,v,w)\in V^{(3)}\subset V^3$ of three distinct
vertices. 

Then either $\Gamma$ or its complementary graph $\Gamma^c$ (having adjacency
matrix $A^c=J-I-A$ where $J$ is the all one matrix and $I$ the identity
matrix) is a disjoint union of at most two complete graphs.
\end{lem}

{\bf Proof.} Up to replacing $\Gamma$ by its complementary graph 
$\Gamma^c$ we can suppose that $\gamma=1$. 
This shows that given any three vertices of $\Gamma$, at least two 
of them are adjacent. 
The graph $\Gamma$ consists thus of at most two connected components.
If a connected component of $\Gamma$ is not a complete graph,
then this component contains two vertices $u,v$ at distance $2$ implying that
$A_{u,v}=0$. Since $u$ and $v$ are at distance $2$ they share a 
common neighbour $w$ for which we have $A_{u,w}=A_{v,w}=1$.
This yields a contradiction since
$A_{u,v}+A_{u,w}+A_{v,w}=2\not\equiv \gamma\pmod 2$.\hfill $\Box$ 

Given a set $E$ we call a function
$\sigma:E^{(2)}\longrightarrow \{\pm 1\}$ {\it symmetric}
if $\sigma(x,y)=\sigma(y,x)$ for all $(x,y)\in E^{(2)}$.  

\begin{prop} \label{corollequiv} Any symmetric function 
$$\sigma:E^{(2)}\longrightarrow 
\{\pm 1\}$$ with 
$$\sigma(a,b)\sigma(b,c)\sigma(a,c)=\gamma\in\{\pm 1\}$$
independent of $(a,b,c)\in E^{(3)}=\{(a,b,c)\in E^3,\ \vert\ 
a\not=b\not= c\not=a\}$
gives rise to a two-partition of $E$.
\end{prop}

{\bf Proof.} Consider the simple graph $\Gamma$ with vertices
$E$ and adjacency matrix having coefficients $A_{x,x}=0$ and
$A_{x,y}=\frac{\sigma(x,y)+1}{2},\ x\not= y$.

The graph $\Gamma$ satisfies the assumptions of Lemma \ref{equivgraphe}
and consists hence, 
up to a sign change of $\sigma$ (which replaces $\Gamma$ by its
complementary graph), of at most two non-empty complete
graphs. The connected components of $\Gamma$ define a two-partition on $E$.
\hfill $\Box$

\begin{rem} \label{equivconstr} 
The two-partition described by Proposition 
\ref{corollequiv}  can
be constructed as follows: set $\gamma=\sigma(a,b)\sigma(b,c)\sigma(a,c)$
for $a\not= b\not= c\not=a$ and choose an element $x_0\in E$.
Up to multiplication by $-1$  the function $\alpha:E\longrightarrow\{\pm 1\}$
defined by $\alpha(x_0)=1$ and $\alpha(x)=\gamma \sigma(x,x_0),\ 
x\not= x_0$ is then independent of the choice of the element $x_0$
and the classes of the associated two-partition are given by
$\alpha^{-1}(1)$ and $\alpha^{-1}(-1)$.
\end{rem}


\section{Symmetric and antisymmetric functions}

A function $\varphi:E^{(l)}\longrightarrow \{\pm 1\}$ (where $E$
is a set)
is {\it $l-$symmetric} or {\it symmetric} if
$$\varphi(\dots,x_{i-1},x_i,x_{i+1},x_{i+2},\dots)=\varphi(x_1,\dots,
x_{i-1},x_{i+1},x_i,x_{i+2},\dots,x_l)$$
for all $1\leq i<l$ and $(x_1,\dots,x_l)\in E^{(l)}$.
We denote by ${\mathcal F}_+(E^{(l)})$ the set of 
all symmetric functions from $E^{(l)}$ to $\{\pm 1\}$.

Similarly, such a  function $\varphi:E^{(l)}\longrightarrow 
\{\pm 1\}$ is {\it $l-$antisymmetric} or {\it antisymmetric} if
$$\varphi(\dots,x_{i-1},x_i,x_{i+1},x_{i+2},\dots)=-\varphi(x_1,\dots,
x_{i-1},x_{i+1},x_i,x_{i+2},\dots,x_l)$$
for all $1\leq i<l$ and $(x_1,\dots,x_l)\in E^{(l)}$.
We denote by ${\mathcal F}_-(E^{(l)})$ the set of 
all antisymmetric functions from $E^{(l)}$ to $\{\pm 1\}$.

The set ${\mathcal F}_\pm(E^{(l)})={\mathcal F}_+(E^{(l)})\cup
{\mathcal F}_-(E^{(l)})$ of all symmetric or
antisymmetric functions from $E^{(l)}$ to $\{\pm 1\}$ 
is a vector space over ${\mathbf F}_2$,
of dimension ${\sharp(E)\choose l}+1$ for $1<l\leq \sharp (E)$ and
$E$ a finite set. The identity element is given by the symmetric constant 
function $E^{(l)}\longrightarrow \{1\}$ and  the 
group-law is the usual product of functions.

We define the {\it signature homomorphisme} $\hbox{sign}:
{\mathcal F}_\pm(E^{(l)})\longrightarrow \{\pm 1\}$ by 
$\hbox{sign}(\varphi)=1$ if
$\varphi\in{\mathcal F}_+(E^{(l)})$ is symmetric and 
$\hbox{sign}(\varphi)=-1$ if
$\varphi\in{\mathcal F}_-(E^{(l)})$ is antisymmetric. The set 
${\mathcal F}_+(E^{(l)})=\hbox{sign}^{-1}(1)$ of all 
symmetric functions on $E^{(l)}$ is of course a subgroup
of ${\mathcal F}_\pm(E^{(l)})$ and the set 
${\mathcal F}_-(E^{(l)})=\hbox{sign}^{-1}(-1)$ of all 
antisymmetric functions on $E^{(l)}$ is a 
free ${\mathcal F}_+(E^{(l)})-$set.


\section{The Orchard morphism}

Given a finite set $E$, the aim of this section is to
construct the {\it Orchard morphism}
$$\rho:{\mathcal F}_\pm (E^{(l)})\longrightarrow {\mathcal E}(E)\ ,$$
a natural group homomorphism
which factors through the quotient group ${\mathcal F}_\pm (E^{(l)})/\pm 1$
where $\pm 1$ denote the obvious constant symmetric functions on $E^{(l)}$. 
Naturality means that $\rho$
is equivariant with respect to the obvious actions of the symmetric group
$\hbox{Sym}(E)$ on ${\mathcal F}_\pm(E^{(l)})$ and ${\mathcal E}(E)$. 

Given a totally ordered set $X$ we denote by ${X\choose k}$ 
the set of all strictly increasing sequences of length $k$ in $X$.

For an arbitrary set $X$, we define ${X\choose k}$ 
by choosing first an arbitrary total order relation on $X$. 

Given a function $\varphi\in {\mathcal F}_\pm 
(E^{(l)})$ where $E$ is finite, we define 
$\sigma_\varphi:E^{(2)}\longrightarrow \{\pm 1\}$ by setting
$$\sigma_\varphi(y,z)=\prod_{(x_1,\dots,x_{l-1})\in{E\setminus \{y,z\}\choose l-1}}
\varphi(x_1,\dots,x_{l-1},y)\ \varphi(x_1,\dots,x_{l-1},z)\ .$$

\begin{prop} \label{orchprop} 
The function $\sigma_\varphi$ is a well-defined symmetric 
function on $E^{(2)}$ such that
$$\sigma_\varphi(a,b)\sigma_\varphi(b,c)\sigma_\varphi(a,c)=
\left(\hbox{sign}(\varphi)\right)^{\sharp(E)-3\choose l-2}$$
for all $(a,b,c)\in E^{(3)}$ where $\hbox{sign}(\varphi)$ is the signature 
homomorphism sending $l-$symmetric functions on $E^{(l)}$ to $1$ and
$l-$antisymmetric functions to $-1$.
\end{prop}

{\bf Proof.} Since every sequence $(x_1,\dots,x_{l-1})\in 
{E\setminus\{y,z\}\choose l-1}$ is involved twice in $\sigma_\varphi(y,z)$, 
the value of $\sigma_\varphi(y,z)$ is 
independent of the choice of a particular total order on $E\setminus\{y,z\}$.
Symmetry ($\sigma_\varphi(y,z)=\sigma_\varphi(z,y)$) of $\sigma_\varphi$ 
is obvious.

Consider now
first an element $(x_1,\dots,x_{l-1})\in {E\setminus\{a,b,c\}\choose 
l-1}$. Such an element contributes always a factor $1$ to the product
$\sigma_\varphi(a,b)\sigma_\varphi(b,c)\sigma_\varphi(a,c)$. The pro\-duct 
$\sigma_\varphi(a,b)\sigma_\varphi(b,c)\sigma_\varphi(a,c)$ 
is thus equal to the the product
over all elements $(x_1,\dots,x_{l-2})\in{E\setminus \{a,b,c\}\choose l-2}$
of factors of the form
$$\begin{array}{l}
\displaystyle \varphi(x_1,\dots,x_{l-2},c,a)\ 
\varphi(x_1,\dots,x_{l-2},c,b)\\  
\displaystyle \varphi(x_1,\dots,x_{l-2},a,b)\ 
\varphi(x_1,\dots,x_{l-2},a,c)\\  
\displaystyle \varphi(x_1,\dots,x_{l-2},b,a)\ 
\varphi(x_1,\dots,x_{l-2},b,c)\end{array}$$
and each of these ${\sharp(E) -3\choose l-2}$ factors yields a 
contribution of $\hbox{sign}(\varphi)$. 
\hfill $\Box$

By Proposition \ref{orchprop} the function $\sigma_\varphi$ satisfies 
the conditions
of Proposition \ref{corollequiv} and gives rise to a two-partition 
$\rho(\varphi)\in {\mathcal E}(E)$. We call the application
$\rho:{\mathcal F}_\pm(E^{(l)})\longrightarrow {\mathcal E}(E)$
defined in this way the {\it Orchard morphism}. Given an
element $\varphi\in{\mathcal F}_\pm (E^{(l)})$ we call 
the two-partition $\rho(\varphi)$ the {\it Orchard-partition}
of $\varphi$. The two-partition $\rho(\varphi)$  gives rise to the 
{\it Orchard-equivalence relation} partitioning the elements of $E$
into at most two {\it Orchard classes}. 

\begin{rem} \label{exotic_homo}
(i) If $l=1$, the two-partition $\rho(\varphi)$
on $E$ given by the
Orchard morphism is the obvious one with classes $\varphi^{-1}(1)$ and
$\varphi^{-1}(-1)$.

\ \ (ii) Consider a $2-$symmetric function $\varphi\in{\mathcal F}_+
(E^{(2)})$ satisfying the condition of Proposition
\ref{corollequiv}. By Proposition \ref{corollequiv} it gives rise 
to a two-partition on $E$. If $E$ is finite, we get a
second two-partition on $E$ by considering the Orchard
morphism $\rho(\varphi)$. An easy computation shows
that these two-partitions coincide if $\sharp (E)$ is odd. If $\sharp( E)$ is
even, the image of the Orchard morphism $\rho(\varphi)\in{\mathcal E}(E)$ 
is trivial for such a function $\varphi$.
\end{rem}

Before stating the main result concerning the Orchard morphism,
we recall the definition of equivariance: Let a group $G$ act
on two sets $X$ and $Y$. An application $\psi:X\longrightarrow
Y$ is {\it $G-$equivariant} if $\psi\left(g(x)\right)=g\left(\psi(x)
\right)$ for all $g\in G$ and $x\in X$. Given a set $E$, the
symmetric group $\hbox{Sym}(E)$ of all bijections of $E$ acts in
an obvious way on the groups ${\mathcal F}_\pm(E^{(l)})$ and
${\mathcal E}(E)$ and it is hence natural to study group homomorphisms
from ${\mathcal F}_\pm(E^{(l)})$ to ${\mathcal E}(E)$ which are
{\it natural}, i.e. $\hbox{Sym}(E)-$equivariant.

\begin{thm} \label{orchthm}
For any finite set $E$ and any natural integer $1\leq l<\sharp(E)$, 
the Orchard morphism
$$\rho:{\mathcal F}_\pm(E^{(l)})\longrightarrow {\mathcal E}(E)$$
is the unique natural group homomorphism which is non-trivial. 
Moreover, $\rho$ factors through the quotient group
${\mathcal F}_\pm^l(E)/\{\pm 1\}$ (where, as always, $\pm 1$ denote the
constant symmetric functions on $E^{(l)}$).
\end{thm}

\begin{rem} \label{remunicity} The Orchard morphism exists and is always 
trivial for $l=\sharp(E)$.

For $n=l=2$ there exists an \lq\lq exotic'' natural homomorphisme
which is non-trivial:
Defining $\rho':{\mathcal F}_\pm(E^{(2)})\longrightarrow {\mathcal E}(E)$
(where $E=\{1,2\}$ has two elements) 
by $\rho'(\varphi)=\left(E=\{1,2\}\right)$ if
$\varphi\in{\mathcal F}_+(E^{(2)})$ and 
$\rho'(\varphi)=\left(E=\{1\}\cup\{2\}\right)$ if
$\varphi\in{\mathcal F}_-(E^{(2)})$ we have a natural homomorphism
distinct from the Orchard morphism (which is trivial in this case).

This failure is due the fact that both
two-partitions on the set $E=\{1,2\}$ are $\hbox{Sym}(E)$-invariant.
However, for finite sets $E$ having more than $2$ elements, only the 
trivial two-partition is $\hbox{Sym}(E)-$invariant. 
\end{rem}

A {\it flip} is a symmetric function
$f_X\in {\mathcal F}_+(E^{(l)})$ such that $f_X^{-1}(-1)\subset E^{(l)}$ 
consists (up to permutation of its elements) of a unique 
sequence $X=(x_1,\dots,x_l)\in{E\choose l}$.
We call the set $X=\{x_1,\dots,x_l\}$ the {\it flipset} of the flip $f_X$.

The set $\{f_X\}_{X\in{E\choose l}}$ of all flips is obviously a 
basis of the subspace ${\mathcal F}_+(E^{(l)})$ of symmetric
functions on $E^{(l)}$.

\begin{lem} \label{fliplem} 
Given a flip $f_X\in {\mathcal F}_+(E^{(l)})$ and an arbitrary
element $\varphi\in{\mathcal F}_\pm(E^{(l)})$ 
we have for $(a,b)\in E^{(2)}$
$$\sigma_\varphi(a,b)\ \sigma_{(\varphi f_X)}(a,b)=-1$$
if and only if exactly one of the elements $a,b$ belongs to $X$.
\end{lem}

{\bf Proof.} In the product defining 
$\sigma_\varphi(a,b)\sigma_{(\varphi f_X)}(a,b)$ the factor
$$\begin{array}{l}
\varphi(x_1,\dots,x_{l-1},a)\varphi(x_1,\dots,x_{l-1},b)\\
\quad (\varphi\ f_X)(x_1,\dots,x_{l-1},a)(\varphi\ f_X)(x_1,\dots,x_{l-1},b)\\
\qquad =f_X (x_1,\dots,x_{l-1},a)\ f_X (x_1,\dots,x_{l-1},b)\end{array}$$
corresponding to $(x_1,\dots,x_{l-1})\in {E\setminus\{a,b\}\choose l-1}$
yields a contribution of $1$ except if $X=\{x_1,\dots,x_{l-1},a\}$ or if 
$X=\{x_1,\dots,x_{l-1},b\}$. This happens at most once and only 
if exactly one of the elements $a,b$ belongs to the set $X$.\hfill$\Box$

Lemma \ref{fliplem} implies easily the following result.

\begin{prop} \label{flipcor}
(i) The classes of the two-partition $\rho(f_X)$ associated to a flip $f_X
\in {\mathcal F}_+(E^{(l)})$ are given by $X$ and $E\setminus X$.

\ \ (ii) If two functions $\varphi,\psi=\varphi\ f_X\in {\mathcal F}_
\pm(E^{(l)})$
differ by a flip then the corresponding equivalence relations
$\rho(\varphi)$ and $\rho(\psi)=\rho(\varphi f_X)$ differ exactly
on the subsets $X\times\left(E\setminus X\right)$ and 
$\left(E\setminus X\right)\times X$ of $E\times E$.
\end{prop}

{\bf Proof of Theorem \ref{orchthm}.} Since the set of all flips generates 
${\mathcal F}_+(E^{(l)})$ and since ${\mathcal F}_-(E^{(l)})$ is a free
${\mathcal F}_+(E^{(l)})-$set, the Orchard morphism $\rho$ behaves well
under composition by assertion (ii) of Proposition \ref{flipcor}. 
Since the equivalence
relation associated to a constant function $\pm 1\in {\mathcal F}_+(E^{(l)})$
is obviously trivial, $\rho$ defines a group homomorphism from the quotient
group ${\mathcal F}_+(E^{(l)})/\{\pm 1\}$ into ${\mathcal E}(E)$.

Equivariance of $\rho$ with respect to $\hbox{Sym}(E)$ is obvious.

We have yet to show that every other natural ($\hbox{Sym}(E)-$equivariant)
homomorphism $\rho':{\mathcal F}_\pm(E^{(l)})\longrightarrow {\mathcal E}(E)$
is either trivial or coincides with the Orchard morphism $\rho$.
 
A flip $f_X$ is clearly invariant under the
subgroup $\hbox{Sym}(X)\times\hbox{Sym}(E\setminus X)\subset \hbox{Sym}(E)$.
If $\sharp(E)>2$, any two-partition invariant under 
$\hbox{Sym}(X)\times\hbox{Sym}(E\setminus X)$ of $E$ is either trivial
or equal to $X\cup(E\setminus X)$. This implies that we have either
$\rho'(f_X)=1$ or $\rho'(f_X)=\rho(f_X)$ for any
$\hbox{Sym}(E)-$equivariant 
homomorphism $\rho':{\mathcal F}_\pm (E^{(l)})\longrightarrow 
{\mathcal E}(E)$. Since $\hbox{Sym}(E)$ acts transitively on the set of
all flips, the first case implies triviality of $\rho'$
restricted to ${\mathcal F}_+(E^{(l)})$ while we have
$\rho'=\rho$ for the restriction onto ${\mathcal F}_+(E^{(l)})$ in
the second case. This conclusion holds also for $\sharp(E)=2$ and $l=1$
as can easily be checked. 

If $\rho'$ restricted to ${\mathcal F}_+(E^{(l)})$ is trivial,
the identity ${\mathcal F}_\pm (E^{(l)})=\varphi\ {\mathcal F}_+(E^{(l)})$
for any $\varphi\in {\mathcal F}_-(E^{(l)})$ shows that $\rho'$
restricted to ${\mathcal F}_-(E^{(l)})$ is constant and hence trivial
for $\sharp(E)>2$ by $\hbox{Sym}(E)-$equivariance. For $\sharp(E)=2$ and 
$l=2$, this conclusion fails as shown by the example of Remark
\ref{remunicity}.

We might hence suppose that $\rho'=\rho$ on 
${\mathcal F}_+(E^{(l)})$. Choose an antisymmetric function
$\varphi\in{\mathcal F}_-(E^{(l)})$. If $n=\sharp(E)$ is odd, choose a cyclic
permutation $\beta\in\hbox{Sym}(E)$ (of maximal length $n$) of $E$ and 
consider
$$\tilde\varphi(x_1,\dots,x_l)=\prod_{j=0}^{n-1}\varphi(\beta^j(x_1),
\beta^j(x_2),\dots,\beta^j(x_l))$$
where $\beta^0(x)=x$ and $\beta^j(x)=\beta(\beta^{j-1}(x))$ for 
$x\in E$. The function $\tilde \varphi:E^{(l)}\longrightarrow \{\pm 1\}$
is antisymmetric on $E^{(l)}$ and invariant under the cyclic subgroup
generated by $\beta\in\hbox{Sym}(E)$. The corresponding
two-partition $\rho'(\tilde \varphi)$ is also 
invariant under the cyclic permutation $\beta$ 
and hence trivial since $n=\sharp(E)$
is odd. The equality ${\mathcal F}_-(E^{(l)})=\tilde \varphi {\mathcal F}_+(
E^{(l)})$ implies now the result.

Suppose now $n=\sharp(E)$ even. Choose an element $z\in E$ and a cyclic 
permutation $\beta$ of all $(n-1)$ elements of $E\setminus\{z\}$.
Setting 
$$\tilde\varphi(x_1,\dots,x_l)=\prod_{j=0}^{n-2}\varphi(\beta^j(x_1),
\beta^j(x_2),\dots,\beta^j(x_l))$$
for a fixed element $\varphi\in{\mathcal F}_-(E^{(l)})$ and reasoning as
above we see that $\rho'(\tilde\varphi)\in {\mathcal E}(E)$ is
either trivial or corresponds to the two-partition
$\{z\}\cup(E\setminus \{z\})$. This implies that the same conclusion 
holds for $\rho'(\tilde\varphi)\rho(\tilde\varphi)$ and the
identity ${\mathcal F}_-(E^{(l)})=\tilde \varphi {\mathcal F}_+(E^{(l)})$
shows that the product $\rho'(\varphi)\rho(\varphi)\in{\mathcal E}(E)$ is
constant for $\varphi\in{\mathcal F}_-(E^{(l)})$. By 
$\hbox{Sym}(E)-$equivariance this is only possible if $n=2$ (cf.
Remark \ref{remunicity}) or if $\rho'(\varphi)\rho(\varphi)$ is trivial
which establishes the Theorem.\hfill $\Box$

\subsection{An easy characterisation of $\rho$ restricted to 
${\mathcal F}_+(E^{(l)})$}

In this subsection we give an elementary description
of $\rho(\varphi)$ for $\varphi\in{\mathcal F}_+(E^{(l)})$ an $l-$symmetric
function.

Given a finite set $E$ and an $l-$symmetric function
$\varphi\in{\mathcal F}_+(E^{(l)})$ we 
consider the function $\mu_\varphi:E\longrightarrow \{\pm 1\}$ defined by
$$\mu_\varphi(x)=\prod_{(x_1,\dots,x_{l-1})\in {E\setminus \{x\}\choose l-1}}
\varphi(x_1,\dots,x_{l-1},x)\ .$$

\begin{prop} The two classes of the Orchard relation 
$\rho(\varphi)$ are given by $\mu_\varphi^{-1}(1)$ and $\mu_\varphi^{-1}(-1)$.
\end{prop}

{\bf Proof.} The result clearly holds for the two $l-$symmetric 
constant functions. The Proposition follows now from the fact that 
$\mu_\varphi$ and $\mu_{\varphi f_X}$ differ exactly on $X$ for a 
flip $f_X\in{\mathcal F}_+(E^{(l)})$.\hfill $\Box$

Another proof can be given by remarking that $\mu_\varphi$ defines a
non-trivial $\hbox{Sym}(E)-$equivariant homomorphism into ${\mathcal E}(E)$
which must be the Orchard homomorphism by unicity.

\begin{rem} Setting 
$$\tilde \mu_\varphi(x)=
\prod_{(x_1,\dots,x_{l})\in {E\setminus \{x\}\choose l}}
\varphi(x_1,\dots,x_{l})\ $$
we have $\tilde \mu_\varphi=\mu_\varphi$, up to a sign given by
$$\prod_{(x_1,\dots,x_{l})\in {E\choose l}}
\varphi(x_1,\dots,x_{l})\ .$$
\end{rem}


\section{Generic configurations of points in ${\mathbf R}^d$}

A finite set ${\mathcal P}=\{P_1,\dots,P_n\}$ of $n$ points in the 
oriented real affine space
${\bf R}^d$ is a {\it generic configuration} if any subset of at most $d+1$
points in ${\mathcal P}$
is affinely independent. Generic configurations of $n\leq d+1$
points in ${\bf R}^d$ are simply vertices of $(n-1)-$dimensional
simplices. For $n\geq d+1$, 
genericity boils down to the fact that any set of $d+1$ 
points in ${\mathcal P}$ spans ${\bf R}^d$ affinely. 

Two generic configurations ${\mathcal P}^1$ and ${\mathcal P}^2$ of
${\mathbf R}^d$ are {\it isomorphic} if there
exists a bijection $\sigma:{\mathcal P}^1\longrightarrow {\mathcal P}^2$
such that all pairs of corresponding $d-$dimensional simplices (with vertices 
$(P_{i_0},\dots,
P_{i_d})\subset {\mathcal P}^1$ and $(\sigma(P_{i_0}),\dots,
\sigma(P_{i_d}))\subset {\mathcal P}^2$) have the same orientations
(given for instance for the first simplex 
by the sign of the determinant of the $d\times d$ matrix with rows 
$P_{i_1}-P_{i_0},\dots,P_{i_d}-P_{i_0}$).

Two generic configurations ${\mathcal P}(-1)$ and ${\mathcal P}(+1)$ 
are {\it isotopic} if there exists a continuous path (with respect to
the obvious topology on ${\mathbf R}^{dn}=\left({\mathbf R}^d\right)^n$)
of generic configurations ${\mathcal P}(t),\ t\in [-1,1],$
which joins them. Isotopic
configurations are of course isomorphic. I ignore to what extend the
converse holds. 

Given a finite generic configuration ${\mathcal P}=\{P_1,\dots,P_n\}\subset
{\mathbf R}^d$ we consider the $(d+1)-$antisymmetric function
$\varphi:{\mathcal P}^{(d+1)}\longrightarrow \{\pm 1\}$ defined by
$$\varphi(P_{i_0},\dots,P_{i_d})=1$$
if
$$\det(P_{i_1}-P_{i_0},P_{i_2}-P_{i_0},\dots,P_{i_d}-P_{i_0})>0$$
and $\varphi(P_{i_0},\dots,P_{i_d})=-1$ otherwise.
The Orchard morphism $\rho(\varphi)\in{\mathcal E}({\mathcal P})$
(extended to be trivial if $n\leq d+1$) provides now
a two-partition of the set ${\mathcal P}$.

The associated equivalence relation can be constructed geometrically as
follows: Given two points $P,Q\in
{\bf R}^d\setminus H$, call an
affine hyperplane $H\subset {\bf R}^d$ {\it separating} 
if $P,Q$ are not in the same connected component of 
${\bf R}^d\setminus H$.
For two points $P,Q$ of a finite generic configuration 
${\mathcal P}=\{P_1,\dots,
P_n\}\subset {\bf R}^d$ we denote by $s(P,Q)$ the number of
separating hyperplanes which are affinely spanned by $d$ distinct elements 
in ${\mathcal P}\setminus \{P,Q\}$. The
number $s(P,Q)$ depends obviously only of the isomorphism type of 
$\mathcal P$ and of $P,Q\in{\mathcal P}$.

\begin{prop} \label{geomdescr}
The equivalence relation $\rho(\varphi)$ on a
finite generic configuration ${\mathcal P}\subset {\mathbf R}^d$
is given by $P\sim Q$ if either
$P=Q$ or if $s(P,Q)\equiv  {n-3\choose d-1}\pmod 2$.
\end{prop}

{\bf Proof.} Given two points $P,Q\in {\mathcal P}$ we have 
$$\sigma(P,Q)=\prod_{(R_1,\dots,R_d)\in{{\mathcal P}\setminus\{P,Q\}
\choose d}}\varphi(R_1,\dots,R_d,P)\varphi(R_1,\dots,R_d,Q)=
(-1)^{\alpha(P,Q)}$$ where $\alpha(P,Q)$ denotes the number of subsets
$(R_1,\dots,R_d)\in{\mathcal P}\setminus\{P,Q\}$ such that the two simplices
with cyclically ordered vertices $(R_1,\dots,R_d,P)$ and $(R_1,\dots,R_d,Q)$ 
have opposite 
orientations. This happens if and only if the affine hyperplane containing
the points
$R_1,\dots,R_d$ separates $P$ from $Q$. We have hence $\alpha(P,Q)=s(P,Q)$
and
$$\sigma_\varphi(P,Q)\sigma_\varphi(Q,R)\sigma_\varphi(P,R)=
(-1)^{\sharp({\mathcal P}) -3 \choose (d+1)-2}$$
and Proposition \ref{geomdescr} follows from  Remark 
\ref{equivconstr}.\hfill $\Box$  

A {\it geometric flip} is a continuous path $$t\longmapsto
{\mathcal P}(t)=(P_1(t),\dots,P_n(t))\in\left({\bf R}^d\right)^n,\ 
t\in [-1,1]$$
with ${\mathcal P}(t)=\{P_1(t),\dots,P_n(t)\}$ generic except for $t=0$
where there exists exactly one subset ${\mathcal F}(0)=
(P_{i_0}(0),\dots,P_{i_d}(0))\subset{\mathcal P}(0)$, 
called the {\it flipset}, of $(d+1)$
points contained in an affine hyperplane spanned by any subset of
$d$ points in
${\mathcal F}(0)$. We require moreover that the simplices
$(P_{i_0}(-1),\dots,P_{i_d}(-1))$ and $(P_{i_0}(1),\dots,P_{i_d}(1))$
carry opposite orientations. Geometrically this means that a point 
$P_{i_j}(t)$ crosses the affine hyperplane spanned by
${\mathcal F}(t)\setminus \{P_{i_j}(t)\}$ for $t=0$.

It is easy to see that two generic configurations ${\mathcal
P}_1,{\mathcal P}_2\subset {\bf R}^d$ having $n$ points
can be related by a continuous path involving at most a finite 
number of geometric flips.

The next result follows directly from the fact that two configurations 
${\mathcal P}(1)$ and ${\mathcal P}(-1)$ related by a geometric flip
give rise to $(d+1)-$antisymmetric functions $\varphi_+,\varphi_-\in
{\mathcal F}_-({\mathcal P}^{(d+1)})$ which differ only by a flip: 

\begin{prop} \label{geomflip} Let ${\mathcal P}(-1),{\mathcal P}(+1)\subset 
{\bf R }^d$ be two generic 
configurations related by a flip with respect to a subset 
${\mathcal F}(t)$ of $(d+1)$ points. 

\ \ (i) If two distinct points $P(t),Q(t)$ are either both contained in 
${\mathcal F}(t)$ or both contained in its complement
${\mathcal P}(t)\setminus {\mathcal F}(t)$ then we have 
$$P(-1)\sim Q(-1)\hbox{ if and only if }P(+1)\sim Q(+1)\ .$$

\ \ (ii) For $P(t)\in {\mathcal F}(t)$ and $Q(t)\not\in {\mathcal
  F}(t)$ we have
 $$P(-1)\sim Q(-1)\hbox{ if and only if }P(+1)\not\sim Q(+1)\ .$$
\end{prop}

\centerline{\epsfysize6cm\epsfbox{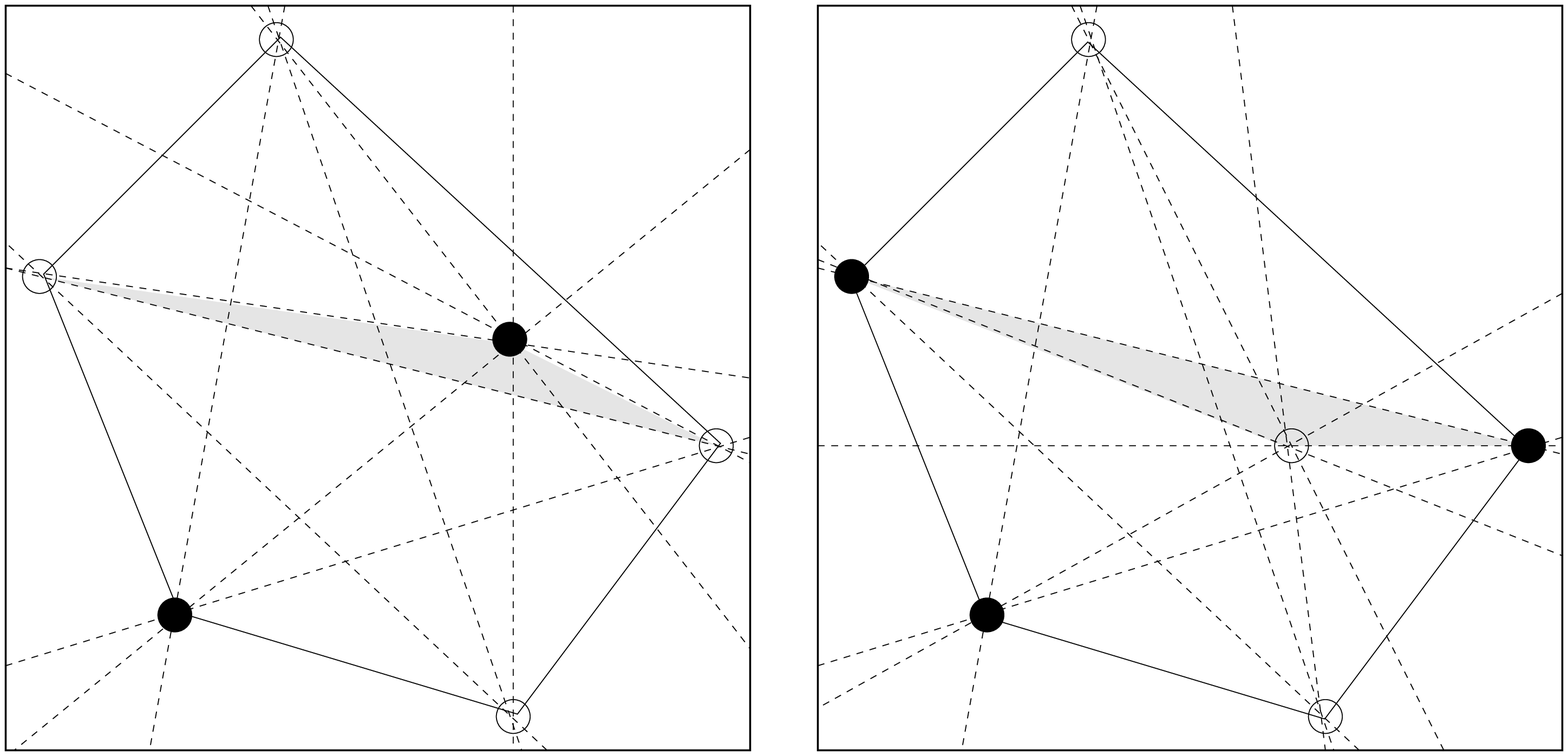}}
\centerline{Figure 2:
Two configurations of 6 points related by a geometric flip} 
\medskip

Proposition \ref{geomflip} suggests also perhaps interesting
problems concerning generic configurations: 
Call two generic
configurations of $n$ points in ${\bf R}^{2d+1}$ {\it
Orchard-equivalent} if they are related by flips
whose flipsets have always exactly $(d+1)$ points in each class.
  
More generally, flips are of different types according to
the number of points of each class involved in the corresponding
flipset. A very special type of flips are the 
{\it monochromatic} ones, defined as involving only vertices of one
class in their flipset. 

Understanding isomorphism classes
of generic configurations
up to flips subject to some restrictions (e.g. only
monochromatic flips or configurations up to
orchard-equivalence in odd dimensions) might be interesting.

We close this section by discussing two further examples.

{\bf Example.} Consider a configuration 
${\mathcal P}\subset {\bf S}^2\subset {\bf R}^3$ consisting of $n$ points 
contained in the Euclidean unit sphere ${\bf S}^2$ and which are generic
as a a subset of ${\bf R}^3$ in the above sense, i.e. $4$ distinct 
points of ${\mathcal P}$ are never contained in a common affine plane
of ${\bf R}^3$. A stereographic projection $\pi: {\bf S}^3\setminus \{N\}
\longrightarrow {\bf R}^2$
with respect to a point $N\in {\bf S}^2\setminus {\mathcal P}$
sends the set ${\mathcal P}\subset {\bf S}^2$ into a set
$\tilde{\mathcal P}=\pi({\mathcal P})\subset {\bf R}^2$ such that $4$ points
of $\tilde {\mathcal P}$ are never contained in a common Euclidean
circle or line of ${\bf R}^2$. The Orchard relation on ${\mathcal P}$ 
can now be seen on $\tilde {\mathcal P}$ as follows: Given two distinct 
points $\tilde P\not=\tilde Q\in\tilde{\mathcal P}$ count the number
$s(\tilde P,\tilde Q)$ of circles or lines determined by $3$ points
in $\tilde {\mathcal P}\setminus\{\tilde P,\tilde Q\}$ which separate them.
Two distinct points $P\not= Q\in{\mathcal P}$ 
are now Orchard-equivalent if and only if
$s(\tilde P,\tilde Q)\equiv {n-3\choose 2}\pmod 2$.
This example can of course be generalised to finite generic
configurations of points on the the $d-$dimensional unit sphere 
${\mathbf S}^d\subset {\mathbf R}^{d+1}$ for $d\geq 2$.

Let ${\mathcal C}$ be a set of continuous real functions on ${\bf R}^k$. 
Suppose ${\mathcal C}$ is a $(d+1)-$dimensional vector
space containing the constant functions. 
Call a set ${\mathcal P}\subset {\bf R}^k$ of $n$ points 
{\it ${\mathcal C}-$generic} if for each subset $S=\{P_{i_1},\dots,
P_{i_d})$ of $d$ distinct points in $\mathcal P$ the set
$$I(S)=\{f\in{\mathcal C}\ \vert\ f(P_{i_j})=0,\ j=0\dots,d\}$$
is a $1-$dimensional affine line 
and all ${n\choose d}$ affine lines in ${\mathcal C}$ of this
form are distinct.

Given $P,Q\in{\mathcal P}$, call a set $S=\{P_{i_1},\dots,
P_{i_d}\}\subset {\mathcal P}\setminus\{P,Q\}$ of $d$ points as above 
{\it ${\mathcal C}-$separating}
(or {\it separating} for short) if $f(P)f(Q)<0$ for any $0\not=f\in I(S)$
and denote by $s_{\mathcal C}(P,Q)$ the number of ${\mathcal C}-$separating
subsets of $\mathcal P$.

\begin{prop} \label{geometric_orchard_relation}
The relation $P\sim_{\mathcal C} Q$ if
either $P=Q$ or
$$s_{\mathcal C}(P,Q)\equiv {n-3\choose d-1}\pmod 2$$
defines an equivalence relation having at most two classes on a 
set ${\mathcal P}=\{P_1,\dots,P_n\}\subset {\bf R}^k$ 
of $n$ points in ${\bf R}^k$ which are ${\mathcal C}-$generic.
\end{prop}

\noindent{\bf Proof.} Consider the linear map $V:{\mathbf R}^k\longrightarrow
{\mathbf R}^d$ defined by $x\longmapsto V(x)=(b_1(x),\dots,b_d(x))$ where
$1,b_1,\dots,b_d$ is a basis of the vector space ${\mathcal C}$.
The image $V({\mathcal P})$ of a ${\mathcal C}-$generic set ${\mathcal P}
\subset {\mathbf R}^k$ is a generic subset of ${\mathbf R}^d$ and
the relation defined by Proposition \ref{geometric_orchard_relation}
coincides with the Orchard-relation described for instance by 
\ref{geomdescr}.\hfill$\Box$

{\bf Examples.} (i) Considering the $(d+1)-$dimensional 
vector space of all affine functions in ${\bf R}^d$, Corollary 1.6
boils down to Theorem 1.1.

\ \ (ii) Consider the $6-$dimensional 
vector space $\mathcal C$ of all polynomial
functions ${\bf R}^2\longrightarrow {\bf R}$ of degree at most 2. A finite subset $
{\mathcal P}\subset {\bf R}^2$ is ${\mathcal C}-$generic if and only if
every subset of five points in $\mathcal P$ defines a unique conic and
all these conics are distinct.

\ \ (iii) Consider the vector space $\mathcal C$ of all polynomials of degree
$< d$ in $x$ together with the polynomials $\lambda y,\ 
\lambda\in{\bf R}$. A subset ${\mathcal P}=\{(x_1,y_1),\dots,(x_n,y_n)\}$ 
with $x_1<x_2<\dots,<x_n$ is ${\mathcal C}-$generic if all ${n\choose 
d}$ interpolation polynomials in $x$ defined by $d$ points of
$\mathcal P$ are distinct.


\section{Orientable sets}

In the following sections we consider a set $E$ together with a fixpoint-free 
involution $\iota:E\longrightarrow E$. We call $\iota$ the 
{\it orientation-reversion}
and the pair $(E,\iota)$ an {\it orientable set}.
The aim of the following sections is to define 
the Orchard morphism for finite orientable sets. In this case,
all groups and homomorphisms are required to be also natural with respect 
to the involution $\iota$. 

Examples of orientable sets 
are for instance antipodal sets of points
in ${\mathbf R}^d\setminus\{0\}$ or points of real Grassmannians endowed with
orientations.

In the sequel we denote by $\pi:E\longrightarrow {\overline E}=E/\iota$
the quotient map $x\longrightarrow {\overline x}=\{x,\iota(x)\}$ onto
the underlying (unoriented) quotient set. The set of all sections
$${\mathcal S}(E,\iota)=
\{s:{\overline E}\longrightarrow E\vert \pi\circ s({\overline x})
={\overline x}\ ,\forall {\overline x}\in {\overline E}\}$$ 
is endowed with a free action of the group $(\pm 1)^{\overline E}$ 
of all functions ${\overline E}\longrightarrow \{\pm 1\}$ if we set
$(\alpha s)({\overline x})=s(x)$ if $\alpha({\overline x})=1$ and 
$(\alpha s)({\overline x})=\iota(s(x))$ otherwise where
$\alpha:{\overline E}\longrightarrow \{\pm 1\}$ and $s\in{\mathcal S}
(E,\iota)$. The quotient 
set ${\mathcal S}(E,\iota)/\{\pm 1\}$ corresponds to orientations
defined up to global reversion (action of $\iota$). We call an element 
of the quotient group ${\mathcal S}(E,\iota)/\{\pm 1\}$ a 
{\it semi-orientation} of the orientable set $(E,\iota)$.

Given an orientable set $(E,\iota)$, its 
automorphism group $\hbox{Sym}(E,\iota)$ is the set of all
$\iota-$equivariant permutations of $E$. Otherwise stated,
a permutation $\pi:E\longrightarrow E$ belongs to $\hbox{Sym}(E,\iota)$
if and only if $\pi(\iota(x))=\iota(\pi(x))$ for all $x\in E$.
As an abstract group, the group $\hbox{Sym}(E,\iota)$ is
easily seen to be isomorphic to the group of all isometries of
the $e-$dimensional regular standard cube $[-1,1]^e\subset {\mathbf R}^e$ 
where $2e=\vert E\vert$ is the cardinality of $E$. This group has $2^e\ e!$
elements and is the wreath product of $\hbox{Sym}({\overline E})$
with $\{\pm 1\}^e$. We have an obvious surjective homomorphisme
$\hbox{Sym}(E,\iota)\longrightarrow \hbox{Sym}({\overline E})$
with kernel $\{\pm 1\}^e$.


\section{Two-sets of orientable sets}

Given an orientable set $(E,\iota)$ it is natural to consider the set $
{\mathcal E}(E,\iota)$ of all two-partitions of $E$
which are invariant under $\iota$. This set
contains the subset ${\mathcal E}(E,\iota_+)$ consisting of all 
two-partitions factoring through $\pi$ and inducing a two-partition
on the quotient set ${\overline E}$. Otherwise stated, two elements
$x$ and $\iota(x)$ in an orbit under $\iota$ belong always 
to the same class.
We call such a two-partition {\it even} since its classes are given 
$\alpha^{-1}(1)$ and $\alpha^{-1}(-1)$ where 
$\alpha:E\longrightarrow\{\pm 1\}$ is an even function with 
respect to the involution $\iota$ (it satisfies
$\alpha(x)=\alpha(\iota(x))$ for all $x\in E$).
Its complement ${\mathcal E}(E,\iota_-)={\mathcal E}(E,\iota)\setminus
{\mathcal E}(E,\iota_+)$, called the {\it odd} two-partitions, 
has equivalence classes defined as 
preimages of an odd function $\alpha:E\longrightarrow\{\pm 1\}$
satisfying $\alpha(x)=-\alpha(\iota(x))$ for all $x$. The set 
${\mathcal E}(E,\iota_-)$ of all odd 
two-partitions on $(E,\iota)$ coincides with the
set of semi-orientations of the orientable set $(E,\iota)$.
Its elements are unordered pairs $\{s,\iota\circ s\}$ 
of complementary sections of the
quotient map $\pi:E\longrightarrow {\overline E}$.

The set
${\mathcal E}(E,\iota)={\mathcal E}(E,\iota_+)\cup {\mathcal E}(E,\iota_-)$ 
obtained by considering all even or odd two-partitions on the 
orientable set $(E,\iota)$ is a vector space (of dimension $\sharp(E)/2$
if $E$ is finite) over ${\mathbf F}_2$. An element of ${\mathcal E}(E,\iota)$
is represented by $\pm \alpha$ where the function
$\alpha:E\longrightarrow \{\pm 1\}$
is either even ($\alpha(\iota x)=\alpha(x)$ for all $x\in E$) or odd 
($\alpha(\iota x)=-\alpha (x)$ for all $x\in E$) 
with respect to $\iota$. The pair $\pm 1$ of constant even functions
represents the identity element and the
group law is the usual product of (pairs of) functions. Given an element
$\{\pm \alpha\}\in {\mathcal E}(E,\iota)$, we define a {\it parity 
homomorphism} ${\mathcal E}(E,\iota)\longrightarrow\{\pm 1\}$ 
by setting $\hbox{parity}(\pm \alpha)=1$ if
$\{\pm \alpha\}\in {\mathcal E}(E,\iota_+)$ is a pair of even functions 
and $\hbox{parity}(\pm \alpha)=-1$ if $\alpha$ is an odd function.

Given an orientable set $(E,\iota)$, we define the set $(E,\iota)^{(l)}$
as the set of all sequences $(x_1,\dots,x_l)\in E^l$ of length $l$
such that $({\overline x_1},\dots,{\overline x_l})\in {\overline E}^{(l)}$.
Otherwise stated, such a sequence satisfies 
$\{x_i,\iota(x_i)\}\not=\{x_j,\iota(x_j)\}$ for $i\not= j$.

\begin{prop} \label{orequivrel}
Any even (respectively odd) symmetric function
$$\sigma:(E,\iota)^{(2)}\longrightarrow \{\pm 1\}$$
with
$$\sigma(a,b)\sigma(b,c)\sigma(a,c)=\gamma\in\{\pm 1\}$$
independent of $(a,b,c)\in
(E,\iota)^{(3)}$ gives rise to an even (respectively odd)
two-partition on $(E,\iota)$.
\end{prop}

{\bf Proof.} Results from Proposition \ref{corollequiv} if $\sigma$ is even.

For $\sigma$ odd, choose a section $s:{\overline E}\longrightarrow E$ and
define the two-partition in the obvious way on the section.
This two-partition extends to a unique odd two-partition on $(E,\iota)$ 
which is independent of the choice of the section $s$. \hfill$\Box$

\begin{rem}
The above equivalence relation can be constructed as follows: Choose 
a fixed base point $x_0\in E$. Set $\alpha(x_0)=1$ and 
$\alpha(\iota(x_0))=\hbox{parity}(\sigma)$ where 
$\hbox{parity}(\sigma)=1$ if $\sigma$ is even and $\hbox{parity}(\sigma)=-1$ 
if $\sigma$ is odd.
For $y\not\in\{x_0,\iota(x_0)\}$ we set $\alpha(y)=\gamma\ \sigma(x_0,y)$
with $\gamma\in\{\pm 1\}$ as in Proposition \ref{orequivrel}.
\end{rem}


\section{Symmetric and antisymmetric functions on orientable sets}

Recall that $(E,\iota)^{(l)}$ denotes the set of all sequences
$(x_1,\dots,x_l)\in E^l$ such that $(\overline{x_1},\dots,\overline{x_l})\in
{\overline E}^{(l)}$.

One defines $l-$symmetric (respectively $l-$antisymmetric) functions on
$(E,\iota)^{(l)}$ in the obvious way as the subset of functions which
are invariant (respectively which change sign) under transposition of two
arguments. 

A symmetric or antisymmetric function $\varphi:(E,\iota)^{(l)}
\longmapsto \{\pm 1\}$ is {\it even} if 
$$\varphi(x_1,x_2,\dots,x_l)=\varphi(\iota (x_1),x_2,\dots,x_l)=
\varphi(x_1,\iota (x_2),x_3,\dots,x_l)=\dots
\ .$$
We denote by ${\mathcal F}_\pm (E,\iota_+)^{(l)}$ the set of all even 
$l-$symmetric or $l-$antisymmetric functions. Notice that there exists
an obvious bijection between ${\mathcal F}_\pm (E,\iota_+)^{(l)}$
and ${\mathcal F}_\pm({\overline E}^{(l)})$.

Such a function is {\it odd} if 
$$\varphi(x_1,x_2,\dots,x_l)=-\varphi(\iota (x_1),x_2,\dots,x_l)=
-\varphi(x_1,\iota (x_2),x_3,\dots,x_l)=\dots\ .$$
The set of all odd symmetric or antisymmetric functions on $(E,\iota)^{(l)}$
will be denoted by ${\mathcal F}_\pm(E,\iota_-)^{(l)}$.

We denote by ${\mathcal F}_\pm(E,\iota)^{(l)}=
{\mathcal F}_\pm (E,\iota_+)^{(l)}\cup {\mathcal F}_\pm (E,\iota_-)^{(l)}$
the set of all even or odd, $l-$symmetric or $l-$antisymmetric 
functions on the orientable set $(E,\iota)$. The set 
${\mathcal F}_\pm(E,\iota)^{(l)}$ is of course a vector space
(of dimension ${\sharp(E)/2\choose l}+2$ is $E$ is finite) over 
${\mathbf F}_2$. The set ${\mathcal F}_\pm(E,\iota_-)^{(l)}$
is a free ${\mathcal F}_\pm(E,\iota_+)^{(l)}-$set.

We define the signature and parity homomorphismes $\hbox{sign},\
\hbox{parity}:{\mathcal F}_\pm
(E,\iota)^{(l)}\longrightarrow \{\pm 1\}$ by 
$$\begin{array}{ll}
\displaystyle \hbox{sign}(\varphi)=1\hbox{ if }
\varphi\in {\mathcal F}_+(E,\iota)^{(l)},\quad&
\hbox{sign}(\varphi)=-1\hbox{ if }
\varphi\in {\mathcal F}_-(E,\iota)^{(l)}\ ,\\
\displaystyle \hbox{parity}(\varphi)=1\hbox{ if }
\varphi\in {\mathcal F}_\pm(E,\iota_+)^{(l)},\quad&
\hbox{parity}(\varphi)=-1\hbox{ if }
\varphi\in {\mathcal F}_\pm(E,\iota_-)^{(l)}\ .\end{array}$$


\section{The Orchard morphism for finite orientable sets}

Given $\varphi\in{\mathcal F}_\pm(E,\iota)^{(l)}$ where $(E,\iota)$ is a finite
orientable set, we define 
$\sigma_\varphi:(E,\iota)^{(2)}\longrightarrow \{\pm 1\}$ 
by setting
$$\sigma_\varphi(y,z)=\prod_{({\overline x_1},\dots,{\overline x_{l-1}})\in{
{\overline E}\setminus \{{\overline y},{\overline z}\}\choose l-1}}
\varphi(x_1,\dots,x_{l-1},y)\ \varphi(x_1,\dots,x_{l-1},z)$$
where $x_1=s({\overline x_1}),\dots,x_{l-1}=s({\overline x_{l-1}})$ 
are obtained using an arbitrary section $s:{\overline E}\longrightarrow E$ of 
the quotient map $\pi:E\longrightarrow {\overline E}=E/\iota$.

\begin{prop} Let $\varphi\in{\mathcal F}_\pm(E,\iota)^{(l)}$ 
be a function and define $\sigma_\varphi$ as above. 

\ \ (i) The function $\sigma_\varphi$ is well defined, symmetric and 
satisfies the identity
$$\sigma_\varphi(a,b)\sigma_\varphi(b,c)\sigma_\varphi(a,c)=
\left(\hbox{sign}(\varphi)\right)^{e-3\choose l-2}$$
for all $(a,b,c)\in (E,\iota)^{(3)}$ where 
$2e=\vert E\vert =2\ \vert {\overline E}\vert$ is the cardinality of $E$.

\ \ (ii) If $\varphi\in{\mathcal F}_\pm(E,\iota_+)^{(l)}$ (i.e. $\varphi$ 
even), then $\sigma_\varphi$ is even.

\ \ (iii) If $\varphi\in{\mathcal F}_\pm(E,\iota_-)^{(l)}$ 
(i.e. $\varphi$ odd), then $\sigma_\varphi$ is even if 
${e-2\choose l-1}\equiv 0\pmod 2$ and odd otherwise.
\end{prop}

{\bf Proof.} Every element ${({\overline x_1},\dots,{\overline x_{l-1}})\in{
{\overline E}\setminus \{{\overline y},{\overline z}\}\choose l-1}}$ is
involved twice in $\sigma_\varphi(y,z)$ thus implying that the final value
is independent of the choosen total order on ${\overline E}\setminus\{y,z\}$
and of the choosen section $s:{\overline E}\longrightarrow E$.

The definition of $\sigma_\varphi$ is obviously symmetric with respect to its 
arguments. 

The proof of the identity 
$\sigma_\varphi(a,b)\sigma_\varphi(b,c)\sigma_\varphi(a,c)=
\left(\hbox{sign}(\varphi)\right)^{
e-3\choose l-2}$ is exactly analogous to the corresponding proof in 
the non-orientable case. 

Assertion (ii) is almost obvious since we have 
$$\varphi(x_1,\dots,x_{l-1},u)=\varphi(x_1,\dots,x_{l-1},\iota(u))$$
for $u\in\{a,b\}$ and $\varphi\in{\mathcal F}_\pm (E^{(l)},\iota_+)$
even.

Assertion (iii) follows from the fact that
$$\varphi(x_1,\dots,x_{l-1},u)=-\varphi(x_1,\dots,x_{l-1},\iota(u))$$
for $u\in\{a,b\}$ and $\varphi\in{\mathcal F}_\pm (E^{(l)},\iota_-)$ odd
and from the observation that the definition of $\sigma_\varphi(a,b)$ involves 
${e-2\choose l-1}$ such factors.\hfill$\Box$

The Orchard morphism 
$\rho:{\mathcal F}_\pm(E^{(l)},\iota)\longrightarrow
{\mathcal E}(E,\iota)$ associates to a function $\varphi\in 
{\mathcal F}_\pm(E^{(l)},\iota)$ the two-partition
in ${\mathcal E}(E,\iota)$ associated to $\sigma_\varphi$ by
Proposition \ref{orequivrel}.

\begin{thm} For $2l<\sharp (E)$ and $\sharp(E)\geq 6$, 
the oriented Orchard morphism is 
the unique $\hbox{Sym}(E,\iota)-$equivariant homomorphism
from ${\mathcal F}_\pm(E,\iota)^{(l)}$ to $
{\mathcal E}(E,\iota)$ which is non-trivial.
\end{thm}

\begin{rem} If $(E,\iota)$ is an orientable set containing $4$
elements $\pm a,\pm b$ (with $\iota$ given by $\iota(a)=-a$ and $
\iota(b)=-b$), there exist several non-trivial natural 
homomorphisms ${\mathcal F}_\pm (E,\iota)^{(l)}\longrightarrow
{\mathcal E}(E,\iota)$ for $l=1,2$.

An example (distinct from the Orchard morphism) for $l=1$ is
given by $\rho'(\varphi)=\hbox{trivial}$ if $\varphi\in
{\mathcal F}(E,\iota_+)^{(1)}$ and $\rho'(\varphi)=
\left(E=\{\pm a\}\cup\{\pm b\}\right)$ if $\varphi\in
{\mathcal F}(E,\iota_-)^{(1)}$.

For $l=2$, one can for instance extend the exotic homomorphism of
the unoriented case (cf. Remark \ref{exotic_homo})
in two ways by choosing an arbitrary even two-partition
as the image $\rho'(\varphi)$ for $\varphi\in{\mathcal F}_+(E,\iota_-)^{(2)}$.
The image $\rho'(\psi)$ for $\psi\in{\mathcal F}_-(E,\iota_-)^{(2)}$
is then the unique remaining two-partition (i.e.
$\rho'(\varphi)\rho'(\psi)=\rho'(\theta)$ with $\theta\in
{\mathcal F}_-(E,\iota_+)^{(2)}$ the unique even non-trivial
two-partition of $(E,\iota)$).
\end{rem}

{\bf Proof.} The proof that $\rho$ defines a homomorphism is as 
in the unoriented case.

The restriction of $\rho$ to the subgroup ${\mathcal F}_\pm
(E,\iota_+)^{(l)}$ consisting only of even functions
coincides with the usual Orchard morphism 
${\mathcal F}_\pm ({\overline E})\longrightarrow {\mathcal E}({\overline
E})$ on ${\overline E}$ and the result holds for this restriction by 
Theorem \ref{orchthm}.

It remains to show unicity of the restriction to 
${\mathcal F}_\pm(E,\iota_-)$ of a natural homomorphisme $\rho'$. The
identity ${\mathcal F}_\pm(E,\iota_-)=\varphi{\mathcal F}_\pm(E,\iota_+)$
 for $\varphi\in {\mathcal F}_\pm(E,\iota_-)$ and 
$\hbox{Sym}(E,\iota)-$equivariance show that such a homomorphism
with trivial restriction to ${\mathcal F}_\pm(E,\iota_+)$ is trivial.

We can thus suppose that $\rho'=\rho$ on ${\mathcal F}_\pm(E,\iota_+)$.
We denote by $e=\sharp({\overline E})$ the halved cardinal of $E$. 

Consider now a section $s:{\overline E}\longrightarrow E$ and the unique
symmetric odd function $\varphi\in{\mathcal F}_+(E,\iota_-)$ defined by
$$\varphi(s({\overline x}_{i_1}),\dots,s({\overline x}_{i_l}))=1$$
for all $({\overline x}_{i_1},\dots,{\overline x}_{i_l})\in
{\overline E}^{(l)}$. $\hbox{Sym}(E,\iota)-$equivariance of $\rho'$
implies that $\rho'(\varphi)$ is either trivial or the semi-orientation
associated to the section $s$. Choose now an element 
${\overline x}\in{\overline E}$ and consider the corresponding function
$\tilde \varphi$ associated as above to the section $\tilde s$ which coincides 
with $s$ on ${\overline E}\setminus\{{\overline x}\}$ and sends
${\overline x}$ to $\iota(s({\overline x}))$. The functions 
$\varphi$ and $\tilde \varphi$ differ by the product of all
${e-1\choose l-1}$ flips with flipsets 
$\{{\overline x},{\overline y}_1,\dots,{\overline y}_{l-1}\}$ where
$( {\overline y}_1,\dots,{\overline y}_{l-1})\in 
{{\overline E}\setminus\{{\overline x}\}\choose l-1}$. 
An element ${\overline y}\in{\overline E}\setminus\{{\overline x}\}$
is involved in ${e-2\choose l-2}$ such flipsets and ${\overline x}$
is involved in ${e-1\choose l-1}={e-2\choose l-2}+{e-2\choose l-1}$
such flipsets. This shows that $\rho'(\varphi)=\rho'(\tilde\varphi)$
if ${e-2\choose l-1}\equiv 0\pmod 2$ and $\hbox{Sym}(E,\iota)-$equivariance
forces $\rho'(\varphi)$ to be even. It coincides thus with the 
Orchard morphism.

If ${e-2\choose l-1}\equiv 1\pmod 2$, the two-partitions 
$\rho'(\varphi)$ and $\rho'(\tilde\varphi)$ differ exactly on
$\pi^{-1}({\overline x})$ and $\hbox{Sym}(E,\iota)-$equivariance
forces $\rho'(\varphi)$ to be the semi-orientation of ${\mathcal E}(E,\iota)$
associated to the section $s$.\hfill $\Box$


\section{Geometric examples}

In this section we discuss a few orientable sets arising from
geometric configurations: finite generic antipodal configurations
of points (or generic configurations of lines through the origin) in
${\mathbf R}^d$ and generic configurations of the real projective 
space ${\mathbf R}P^d$. 

A finite {\it antipodal set} of ${\mathbf R}^d$ is a finite subset 
${\mathcal P}\subset {\mathbf R}^d\setminus\{0\}$ invariant
under the involution $x\longmapsto\iota(x)=-x$.  We call such a set 
{\it generic} if the linear span of any subset $\{\pm x_1,\dots,\pm x_k\}
\subset {\mathcal P}$ is of dimension $k$ for $k\leq d$. 
We get an element
$\varphi\in{\mathcal F}_-(E,\iota)^{(d)}$ by considering the 
$\hbox{sign}\in\{\pm 1\}$ of $$\det(x_1,\dots,x_d)$$
for $(x_1,\dots,x_d)\in({\mathcal P},\iota)^{(d)}$ (where 
$\det(x_1,\dots,x_d)$ denotes the non-zero determinant of the $d\times
d$ matrix with rows $x_1,\dots,x_d$).

Applying the oriented Orchard morphism $\rho$ of the preceeding
section to $\varphi$ we get a two-partition $\rho(\varphi)\in
{\mathcal E}({\mathcal P},\iota)$. Obviously, $\rho(\varphi)$ 
remains the same by rescaling each pair $\pm x\in{\mathcal P}$ 
by some strictly positive constant $\lambda_x\in{\mathbf R}_{>0}$.

We may rescale such an antipodal set in order to lie on the 
Euclidean sphere ${\mathbf S}^{d-1}=\{x\in{\mathbf R}^d\ \vert\
\ \parallel x\parallel =1\}\subset {\mathbf R}^d$.
Similarly, we can interprete ${\mathcal P}$ as a set ${\mathcal L}$
of lines (defined by opposite pairs $\pm x\in{\mathcal P}$). The
Orchard morphism $\rho(\varphi)$ endows then such a generic finite set
of lines either with a two-partition (in the case where
${\sharp({\mathcal L})-2\choose d-1}\equiv 0\pmod 2$) or with a 
semi-orientation (if  ${\sharp({\mathcal L})-2\choose d-1}\equiv 1\pmod 2$).

A finite subset ${\mathcal P}\subset {\mathbf R}P^d$ of the real 
projective space is generic if its completed preimage 
${\mathcal L}={\overline{\pi^{-1}({\mathcal P})}}\subset {\mathbf R}^{d+1}$ is 
a finite set of generic lines in ${\mathbf R}^{d+1}$. If 
${\sharp({\mathcal P})-2\choose d}\equiv 0\pmod 2$ we get a 
two-partition on such a set $\mathcal P$ by applying the Orchard morphism to 
${\mathcal L}$.

In the case where the Orchard morphism endows ${\mathcal L}$ with a 
semi-orientation, we get also an interesting structure on ${\mathcal P}$ as
follows:

Any pair $P,Q\in{\mathcal P}$ of distinct points defines two connected 
components on $L_{P,Q}\setminus\{P,Q\}$ where $L_{P,Q}\subset {\mathbf R}P^d$
denotes the projective line containing $P$ and $Q$. One of these 
connected components is now
selected by a semi-orientation on ${\mathcal L}$ by choosing
the connected component of $L_{P,Q}\setminus\{P,Q\}$ 
whose preimage in 
${\mathbf S}^d\subset {\mathbf R}^{d+1}=\pi^{-1}({\mathbf R}P^d)\cup \{0\}$ 
joins elements of $\pi^{-1}({
\mathcal P})$ which are in the same class. We get in this way an
immersion of the complete graph $K_{\mathcal P}$ with vertices
${\mathcal P}$ into the projective space ${\mathcal R}P^d$. It is
straightforward to show that this immersion is homologically trivial:
each cycle of $K_{\mathcal P}$ is immerged in a contractible way
into ${\mathbf R}P^d$. 


\begin{rem} A preliminary version of this paper (cf. \cite{B}) contained also 
a section concerning simple arrangements of (pseudo)lines in
the projective plane. 
The corresponding invariants (two-partitions and semi-orientations)
are however not based on the Orchard-morphism but use only
Proposition \ref{orequivrel}. 
They are thus not really related to the topic of this
text and will be discussed elsewhere.
\end{rem}

I would like to thank many people who where interested in this work,
especially M. Brion, P. Cameron, 
E. Ferrand, P. de la Harpe and A. Marin for their remarks and comments.

\vskip1cm

Roland Bacher

INSTITUT FOURIER

Laboratoire de Math\'ematiques

UMR 5582 (UJF-CNRS)

BP 74

38402 St MARTIN D'H\`ERES Cedex (France)
 
e-mail: Roland.Bacher@ujf-grenoble.fr

\end{document}